\documentclass[12pt,a4paper]{article}
\usepackage{amsmath}
\usepackage{amssymb}
\usepackage{t1enc}
\usepackage[latin1]{inputenc}
\usepackage[english]{babel}
\usepackage{amsfonts}
\usepackage{latexsym}
\usepackage{bm}
\usepackage{setspace}
\usepackage[toc,page]{appendix}
\usepackage{graphicx}
\usepackage{enumerate}
\usepackage[numbers,sort&compress]{natbib}
\usepackage{tikz}

\usepackage{theoremref}
\usepackage{lineno}

\usepackage{mathtools}

\newtheorem{theorem}{Theorem}%[section]

\newtheorem{lemma}{Lemma}

\newtheorem{const}{Construction}

\allowdisplaybreaks

\begin{document}

\title{Isolation of $k$-cliques in uniform hypergraphs}
%\vspace{2cm}
\author{\\[5mm]
Mattea Baluci\\[5mm]
Department of Mathematics \\
Faculty of Science \\
University of Malta\\
Malta\\
\normalsize{\texttt{mattea.baluci.21@um.edu.mt}}
\and 
\\[5mm]
Peter Borg \\ [5mm]
Department of Mathematics \\
Faculty of Science \\
University of Malta\\
Malta\\
{\normalsize \texttt{peter.borg@um.edu.mt}}
}

\date{}%{\today} 
\maketitle

\begin{abstract}
\noindent A copy of a hypergraph $F$ is called an $F$-copy. For any integers $k$ and $r$ with $2 \leq r \leq k$, let $K_k^r$ denote the complete $r$-uniform hypergraph whose vertex set is $[k] = \{1, 2, \dots, k\}$ (that is, the edges of $K_k^r$ are the $r$-element subsets of $[k]$). Given an $r$-uniform $m$-edge hypergraph $H$, the $K_k^r$-isolation number of $H$, denoted by $\iota(H, K_k^r)$, is the size of a smallest subset $D$ of the vertex set of $H$ such that the closed neighbourhood $N[D]$ of $D$ intersects the vertex sets of the $K_k^r$-copies contained by $H$ (equivalently, $H-N[D]$ contains no $K_k^r$-copy). We prove that if $H$ is connected and not a $K_k^r$-copy, then 
$$\iota(H, K^r_k) \le \left \lfloor \frac{(r-1)m + 1}{(r-1)(\binom{k}{r}+1) + 1} \right \rfloor.$$ 
The bound is best possible. This is a solution to a problem of Li, Zhang and Ye. We also provide an alternative proof for $r < k$. Both proofs introduce new ideas for establishing bounds on isolation parameters.

\end{abstract}
\section{Introduction}

The graph isolation problem was introduced by Caro and Hansberg in \cite{CaHa15, CaHa17}. It is a natural and appealing generalization of the classical graph domination problem \cite{C, CH, HHS, HHS2, HL, HL2}. Similarly to the historical development of domination, isolation has rapidly become a very active field of investigation and is expanding in various directions. The seminal paper \cite{CaHa17} addressed this problem from various angles and posed many pertinent problems, some of which are solved in \cite{Borgcycles, BFK, Borgrsc}. It has motivated research at an increasing rate. Many of the known isolation results are referenced in \cite{Borgrsc, Borgisdom}. The more general hypergraph setting was recently treated by Li, Zhang and Ye in \cite{LZY}, in which they posed a problem concerning isolation of $k$-cliques in $r$-uniform hypergraphs. In this paper, we provide a solution, introducing new ideas and establishing methods of dealing with the complexity of the hypergraph setting. We believe that this paves the way for further advances in the study of isolation. 

Unless stated otherwise, we use capital letters such as $X$ to denote sets or graphs, and small letters such as $x$ to denote non-negative integers or elements of a set. The set of positive integers is denoted by $\mathbb{N}$.  For $n \in \{0\} \cup \mathbb{N}$, $[n]$ denotes the set $\{i \in \mathbb{N} \colon i \leq n\}$. Note that $[0]$ is the empty set $\emptyset$. Arbitrary sets are taken to be finite. A set of sets is called a \emph{family}. A set of size $k$ is called a \emph{$k$-element set} or simply a \emph{$k$-set}. For a set $X$, the \emph{power set of} $X$ (the family of subsets of $X$) is denoted by $2^{X}$, and the family of $k$-element subsets of $X$ is denoted by $\binom{X}{k}$ (that is, $\binom{X}{k} = \{A \subseteq X \colon |A| = k\}$). For standard terminology in graph theory, we refer the reader to \cite{BOL, West}; most of the notation and terminology used here is defined in \cite{Borgcycles}.

A \emph{hypergraph} $H$ is a pair $(X, Y)$ such that $X$ is a set denoted by $V(H)$ and called the \emph{vertex set of $H$}, and Y is a subfamily of $2^{X}$ denoted by $E(H)$ and called the \emph{edge set of $H$}. A member of $V(H)$ is called a \emph{vertex of $H$}, and a member of $E(H)$ is called a \emph{hyperedge of $H$} or simply an \emph{edge of $H$}. If $|V(H)| = n$, then $H$ is called an \emph{$n$-vertex hypergraph}. If $|E(H)| = m$, then $H$ is called an \emph{$m$-edge hypergraph}. If $v,w \in e \in E(H)$ and $v \neq w$, then $w$ is called a \emph{neighbour of v in $H$}. If $v \in e \in E(H)$, then $e$ is said to be \emph{incident to v in $H$}. For $v \in V(H)$, $N_H(v)$ denotes the set of neighbours of $v$ in $H$, $N_H[v]$ denotes the closed neighbourhood $\{v\} \cup N_H(v)$ of $v$, and $E_H(v)$ denotes the set of edges of $H$ that are incident to $v$. If $F$ is a hypergraph such that $V(F) \subseteq V(H)$ and $E(F) \subseteq E(H)$, then $F$ is called a \emph{subhypergraph of} $H$, and we say that $H$ \emph{contains} $F$. For $X \subseteq V(H)$, $N_H[X]$ denotes the closed neighbourhood $\bigcup_{v \in X} N_H[v]$ of $X$, $H[X]$ denotes the subhypergraph of $H$ induced by $X$ (that is, $H[X] = (X, E(H)\cap 2^{X})$), and $H - X$ denotes the hypergraph obtained from $H$ by deleting the vertices in $X$ (that is, $H-X = H[V(H)\setminus X]$). With a slight abuse of notation, for $Y \subseteq E(H)$, let $H-Y$ denote the hypergraph obtained from $H$ by removing all edges in $Y$ (that is, $H-Y= (V(H), E(H)\setminus Y)$). Where no confusion arises, the subscript $H$ may be omitted from any notation that uses it; for example, $N_H(v)$ may be abbreviated to $N(v)$. If $E(H) \subseteq \binom{V(H)}{r}$, then $H$ is said to be \emph{$r$-uniform}. A \emph{graph} is a $2$-uniform hypergraph. An $r$-uniform hypergraph is also called an \emph{$r$-graph}.

If $H_1, \dots, H_t$ are hypergraphs whose vertex sets are pairwise disjoint, then $H_1, \dots, H_t$ are said to be \emph{vertex-disjoint}. A hypergraph $H$ is said to be \emph{connected} if for every $v,w \in V(H)$ with $v \neq w$, there exist some $e_1, \dots, e_t \in E(H)$ such that $v \in e_1,\space w \in e_t$ and $e_i \space \cap \space e_{i+1}\space \neq \space \emptyset$ for each $i \in [t-1]$. A \emph{component of $H$} is a maximal connected subhypergraph of $H$ (that is, no other connected subhypergraph of $H$ contains it). Clearly, the components of $H$ are vertex-disjoint, and their union is $H$.

If $F$ and $H$ are hypergraphs, $f : V(F) \rightarrow V(H)$ is a bijection and $E(H) = \{\{f(v) : v \in e\} \space\colon e \in E(F)\}$, then we say that $H$ is a \emph{copy of $F$} or that $H$ is \emph{isomorphic to $F$}, and we write $H \simeq F$.  Thus, a copy of $F$ is obtained by relabeling the vertices of $F$. A copy of $F$ is also called an \emph{$F$-copy}. The \emph{complete $r$-graph} $([k], \binom{[k]}{r})$ is denoted by $K^r_k$ and called a \emph{$k$-clique}. We abbreviate $K^2_k$ to $K_k$. 

If $D \subseteq  V(H)=N[D]$, then $D$ is called \emph{a dominating set of} $H$. The size of a smallest dominating set is called the \emph{domination number of} $H$ and denoted by $\gamma(H)$. If $\mathcal{F}$ is a set of hypergraphs and $F$ is a copy of a hypergraph in $\mathcal{F}$, then we call $F$ an \emph{$\mathcal{F}$-graph}. If $D \subseteq V(H)$ such that $N[D]$ intersects the vertex sets of the $\mathcal{F}$-graphs contained by $H$, then $D$ is called an \emph{$\mathcal{F}$-isolating set of $H$}. Note that $D$ is an $\mathcal{F}$-isolating set of $H$ if and only if $H-N[D]$ contains no $\mathcal{F}$-graph. The size of a smallest $\mathcal{F}$-isolating set of $H$ is denoted by $\iota(H,\mathcal{F})$. If $\mathcal{F}=\{F\}$, then we may replace $\mathcal{F}$ in these defined terms and notation by $F$.  Clearly, for $r \geq 2$, $D$ is a $K_1^r$-isolating set of $H$ if and only if $D$ is a dominating set of $H$, so $\gamma(H) = \iota(H, K_1^r)$. Trivially, $\iota(H, \mathcal{F}) \leq \gamma(H)$. 
 
As mentioned above, the study of isolating sets was introduced by Caro and Hansberg about a decade ago in \cite{CaHa15, CaHa17}. Consider any connected $n$-vertex $m$-edge graph $G$. One of the earliest domination results is the inequality $\gamma(G) \leq n/2$ of Ore \cite{Ore} for $n \ge 2$ (see \cite{HHS}). While deleting the closed neighbourhood of a dominating set of $G$ yields a graph with no vertices, deleting the closed neighbourhood of a $K_2$-isolating set of $G$ yields a graph with no edges. A $K_2$-isolating set is also called a \emph{vertex-edge dominating set} \cite{BCHH,Lewis,LHHF,Peters,TJ,Z}. Caro and Hansberg \cite{CaHa17} proved that $\iota(G,K_2) \leq n/3$ unless $G$ is a $2$-clique or a $5$-cycle. This solved a problem in \cite{BCHH} and was independently proved by \.{Z}yli\'{n}ski \cite{Z}. Solving a problem of Caro and Hansberg \cite{CaHa17}, Borg \cite{Borgcycles} proved that if $\mathcal{F}$ is the set of cycles and $G$ is not a $3$-clique, then $\iota(G,\mathcal{F}) \leq n/4$ (\cite{Borgcon} provides a stronger version). A special case of this result is that $\iota(G, K_3) \leq n/4$. Another problem posed in \cite{CaHa17} essentially was to determine the best possible upper bound on $\iota(G,K_k)$ in terms of $n$ for any $k \geq 1$. This was solved by Borg, Fenech and Kaemawichanurat \cite{BFK}, who proved that 
\begin{equation} \iota(G,K_k) \le \frac{n}{k+1} \label{BFK_bound_n}
\end{equation}
unless $G$ is a $k$-clique or $k = 2$ and $G$ is a $5$-cycle. The above-mentioned results of Ore, Caro and Hansberg, and Borg are the cases $k = 1, 2, 3$, respectively. The authors of \cite{BFK} subsequently proved in \cite{BFK2} that if $G$ is not a $k$-clique, then 
\begin{equation} \iota(G,K_k) \leq \frac{m+1}{\binom{k}{2} + 2}. \label{BFK_bound_m}
\end{equation}
All the bounds above are sharp. Other isolation bounds of this kind in terms of $n$ or $m$ are given in \cite{BBS, Borgcon, Borgrsc, Borgrc2, Borgisdom, CQ, CX, CZ, Yan, ZW4}. It is worth mentioning that domination and isolation have been particularly investigated for maximal outerplanar graphs \cite{BK, BK2, CaWa13, CaHa17, Ch75, DoHaJo16, DoHaJo17, HeKa18, LeZuZy17, Li16, MaTa96, KaJi, To13}, mostly due to connections with Chv\'{a}tal's Art Gallery Theorem \cite{Ch75}.

Not much is known about domination or isolation of $r$-graphs for $r \geq 3$. Consider any connected $n$-vertex $m$-edge $r$-graph $H$. Recall that $\gamma(H) = \iota(H, K_1^r)$. For $r \in \{3, 4\}$, it is shown in \cite{BHT, HS} that $\gamma(H) \leq n/r$, and that this bound is sharp. For $r = 5$, it is shown in \cite{BHT} that $\gamma(H) \leq 2n/9$. Recently, Li, Zhang and Ye \cite{LZY} provided an initial and rather generic study of isolating sets of $r$-graphs, and they posed the problem of establishing bounds on $\iota(H, K_k^r)$. Borg \cite{Borghyp} showed that (\ref{BFK_bound_n}) yields the more general result that for $2 \leq r \leq k$,
\begin{equation} \iota(H,K_k^r) \le \frac{n}{k+1} \label{Bhyp_bound_n}
\end{equation}
unless $H$ is a $k$-clique or $k = 2$ and $G$ is a $5$-cycle. This bound is sharp. In this paper, we provide another solution to the problem of Li, Zhang and Ye by establishing a best possible bound on $\iota(H, K_k^r)$ in terms of $m$, generalizing (\ref{BFK_bound_m}). In order to state our result in full, we first construct an $r$-graph that attains the bound.

\begin{const} \label{Construction1} \emph{Let $m,k,r \in \mathbb{N}$ with $2 \leq r \leq k$. %and $(m, k, r) \neq (1, r, r)$. 
Let $t = (r-1)(\binom{k}{r} + 1) + 1$. By the division algorithm, $(r-1)m + 1 = qt + s$ for some non-negative integers $q$ and $s$ with $s \leq t-1$. If $q = 0$, then let $B_{m,k,r}$ be a connected $m$-edge $r$-graph. Suppose $q \geq 1$. Let $F_1, \dots, F_q$ be pairwise vertex-disjoint copies of $K^r_k$, and let $e_1, \dots, e_q$ be $r$-element sets such that $V(F_1) \cup e_1, \dots, V(F_q) \cup e_q$ are pairwise disjoint and $|e_i \cap V(F_i)| = 1$ for each $i \in [q]$. For each $i \in [q]$, let $v_i \in e_i \setminus V(F_i)$, and let $w_i$ be the member of $e_i \cap V(F_i)$. Let $p = \left \lfloor \frac{q-1}{r-1} \right \rfloor$ and $p' = \left \lfloor \frac{s}{r-1} \right \rfloor$. For each $j \in [p]$, let $R_j = \{v_i \colon 1 + (j-1)(r-1) \leq i \leq 1 + j(r-1)\}$. If $p \neq \frac{q-1}{r-1}$ (so $1 + (p-1)(r-1) < q$), then let $R_{p+1} = \{v_i \colon 1 + (p-1)(r-1) \leq i \leq q\} \cup R'$ such that $R' \subseteq e_q \setminus \{v_q, w_q\}$ and $|R_{p+1}| = r$. If $p = \frac{q-1}{r-1}$, then let $R_{p+1} = e_{q}$. For each $i \in [q]$, let $H_i$ be the $r$-graph with $V(H_i) = V(F_i) \cup e_i$ and $E(H_i) = E(F_i) \cup \{e_i\}$. Let $T = (\bigcup_{j=1}^{p+1} R_j, \{R_j \colon j \in [p+1]\})$. Let $T'$ be an $r$-graph such that $V(T') \cap \bigcup_{i=1}^{q} V(H_i) = e_q \setminus \{w_q\}$ and $|E(T')| = p'$. Let $B_{m,k,r}$ be the connected $r$-graph $\left( V(T') \cup \bigcup_{i=1}^q V(H_i), E(T) \cup E(T') \cup \bigcup_{i=1}^q E(H_i) \right)$. Note that $|E(B_{m,k,r})| = q(\binom{k}{r} + 1) + \left \lceil \frac{q-1}{r-1} \right \rceil + \left \lfloor \frac{s}{r-1} \right \rfloor$. Since $(r-1)m + 1 = qt + s$, $m = q(\binom{k}{r} + 1) + \frac{q + s - 1}{r-1}$, so $r-1$ divides $q + s - 1$, and hence $\frac{q + s - 1}{r-1} = \left \lceil \frac{q-1}{r-1} \right \rceil + \left \lfloor \frac{s}{r-1} \right \rfloor$. Thus, $|E(B_{m,k,r})| = m$.}
\end{const}

We can now state our result, which is proved in the next section. 

\begin{theorem} \label{edgethm}
If $2 \leq r \leq k$, $H$ is a connected $m$-edge $r$-graph, and $H \not\simeq K_k^r$, then
\[\iota(H,K^r_k) \le  \left \lfloor \frac{(r-1)m + 1}{(r-1)(\binom{k}{r}+1) + 1} \right \rfloor .
\]
Moreover, the bound is attained if $H = B_{m,k,r}$.
\end{theorem}
An alternative proof for $r < k$ is also provided in the next section. Both proofs introduce new ideas for establishing bounds on isolation parameters. This provides a key to the solution of similar problems in isolation theory. %hypergraph isolation problems and hence contribute to the rapid expansion of isolation theory. 

\section{Proofs}
We now start working towards proving Theorem~\ref{edgethm}. 

For any hypergraph $H$, let $\mathcal{C}_k^r(H)$ denote the set of $K_k^r$-copies contained by $H$, let ${\rm C}(H)$ denote the set of components of $H$, and let ${\rm C}_k^r(H)$ denote the set of components of $H$ that are $K_k^r$-copies (that is, ${\rm C}_k^r(H) = \mathcal{C}_k^r(H) \cap {\rm C}(H)$). % $\{I \in {\rm C}(H) \colon I \simeq K_k^r\}$.
    
\begin{lemma} \label{lemma1} If $2 \leq r \leq k$, $H$ is an $r$-graph, $X \subseteq V(H)$, $Y \subseteq N[X]$, $Z \subseteq \{e \in E(H) \colon e \cap N[X] \neq \emptyset\}$, and $S$ is a $K_k^r$-isolating set of $(H-Z)-Y$, then $X \cup S$ is a $K_k^r$-isolating set of $H$, and
\begin{equation} \iota(H, K^k_r) \leq |X| + \iota((H - Z) - Y, K^r_k). \nonumber
\end{equation}
\end{lemma}
\textbf{Proof.} We may assume that $|S| = \iota((H-Z)-Y)$. For each $C \in \mathcal{C}_k^r(H) \setminus \mathcal{C}_k^r((H - Z) - Y)$, we have $V(C) \cap Y \neq \emptyset$ or $E(C) \cap Z \neq \emptyset$, so $V(C) \cap N[X] \neq \emptyset$.  Thus, $X \cup S$ is a $K_k^r$-isolating set of $H$. The result follows.~\hfill{$\Box$} 

\begin{lemma} \label{lemma2} If $2 \leq r \leq k$, $H$ is an $r$-graph, %For $2 \leq r \leq k$ and any $r$-graph $H$, 
and for each $I \in {\rm C}(H)$, $D_I$ is a $K_k^r$-isolating set of $I$, then $\bigcup_{I \in {\rm C}(H)} D_I$ is a $K_k^r$-isolating set of $H$, and
$$\iota(H, K^r_k) = \sum_{I \in {\rm C}(H)} \iota(I, K^r_k ).$$
\end{lemma}
\textbf{Proof.} We may assume that $|D_I| = \iota(I, K_k^r)$ for each $I \in {\rm C}(H)$. Since $r \leq k$, $K_k^r$ is connected. Consequently, for each $C \in \mathcal{C}_k^r(H)$, $V(C)$ does not intersect more than one component of $H$, and hence $V(C) \subseteq V(J)$ for some $J \in {\rm C}(H)$. Thus, $\bigcup_{I \in {\rm C}(H)} D_I$ is a $K_k^r$-isolating set of $H$, and hence $\iota(H, K^r_k) \le \sum_{I \in {\rm C}(H)} \iota(I, K^r_k)$. Let $D$ be a smallest $K_k^r$-isolating set of $H$. For each $I \in {\rm C}(H)$, $D \cap V(I)$ is a $K_k^r$-isolating set of $I$. We have $\sum_{I \in {\rm C}(H)} \iota(I, K^r_k) \le \sum_{I \in {\rm C}(H)} |D \cap V(I)| = |D| = \iota(H, K^r_k)$. The result follows.~\hfill{$\Box$}
\\

\noindent\textbf{Proof 1 of Theorem~\ref{edgethm}.} Let $k_r = \binom{k}{r}$ and $t = (r - 1)(k_r + 1) + 1$. We have $|E(K_k^r)| = k_r$.

We first settle the second part of the theorem. Thus, suppose that $H$ is the $r$-graph $B_{m,k,r}$ in Construction~\ref{Construction1}. If $m < k_r$, then $H$ contains no $K_k^r$-copy, so $\iota(H,K_k^r) = 0 = \lfloor (m(r-1) + 1)/t \rfloor$. If $m = k_r$, then the same holds since $H \not \simeq K_k^r$. Suppose $m \geq k_r + 1$. Then, $q \geq 1$. For each $v \in V(H)$, $N[v]$ does not intersect more than one of the vertex sets of the $k$-cliques $F_1, \dots, F_q$, so $\iota(H, K_k^r) \geq q$. Since $\{w_1, \dots, w_q\}$ is a $K_k^r$-isolating set of $H$, $\iota(H, K_k^r) = q$. Now $q$ is the bound in the theorem.

Using induction on $m$, we now prove that the bound in the theorem holds. The bound is trivial if $m \leq 2$. Suppose $m \ge 3$. If $H$ has no $k$-cliques, then $\iota(H,K^r_k) = 0 < ((r-1)m + 1)/t$. Suppose that $H$ has a $k$-clique $C$. Since $H$ is connected and $H$ is not a $k$-clique, $N[v] \setminus V(C) \neq \emptyset$ for some $v \in V(C)$. Let $u \in N[v] \setminus V(C)$, and let $e_u \in E(H)$ such that $v, u \in e_u$. Let $A_v = E(H[N[v]])$. Since $v \in V(C)$, we have $V(C) \subseteq N[v]$, so $E(C) \subseteq A_v$. Since $u \notin V(C)$, $e_u \in A_v \setminus E(C)$. Thus,
\begin{equation} |A_v| \ge |E(C) \cup \{e_u\}| = k_r + 1. \label{eq1}     
\end{equation}
If $V(H) = N[v]$, then $\{v\}$ is a $k$-clique isolating set of $H$, so $\iota(H,K^r_k) = 1 \leq ((r-1)m + 1)/t$ by (\ref{eq1}).
Suppose $V(H)\neq N[v]$. Let $H' = H - N[v]$. Then, $V(H') \neq \emptyset$. Let $\mathcal{H} = {\rm C}(H')$. For each $Q \in \mathcal{H}$, let $A_Q = \{e \in E(H) \colon e \cap N(v) \neq \emptyset \neq e \cap V(Q)\}$. 

Suppose $e \notin A_v \cup \bigcup_{Q \in \mathcal{H}} (E(Q) \cup A_Q)$ for some $e \in E(H)$. Let $\mathcal{H}_e = \{Q \in \mathcal{H} \colon e \cap V(Q) \neq \emptyset\}$. Since $e \notin A_v$, we have $e \setminus N[v] \neq \emptyset$, so $ \mathcal{H}_e \neq \emptyset$. Let $Q_e \in \mathcal{H}_e$. Suppose $v \in e$. Then, $e \subseteq N[v]$. Thus, we have $e \in A_v$, a contradiction. Therefore, $v \notin e$. Since $e \notin A_{Q_e}$, $e \cap N(v) = \emptyset$. Thus, $e \subseteq V(H')$. Let $H'_e = (\bigcup_{Q \in \mathcal{H}_e} V(Q), \{e\} \cup \bigcup_{Q \in \mathcal{H}_e} E(Q))$. Then, $H'_e$ is a connected subhypergraph of $H'$. Since $e \notin E(Q_e)$, we have $E(Q_e) \subsetneq E(H'_e)$, which contradicts $Q_e \in \mathcal{H}$. Therefore,
\begin{equation} E(H) = A_v \cup \bigcup_{Q \in \mathcal{H}} (E(Q) \cup A_Q). \label{H_edges_1}
\end{equation}

For any $Q \in \mathcal{H}$ and $x \in N(v)$, we say that $Q$ is \emph{linked to $x$} if $x \in e$ for some $e \in A_Q$, and we say that $e$ \emph{links $Q$ to $x$} and that $e$ \emph{links $Q$ to $N(v)$}. Since $H$ is connected, each member of $\mathcal{H}$ is linked to at least one member of $N(v)$. Thus, for each $Q \in \mathcal{H}$, we have $x_Q, y_Q \in e_Q$ for some $x_Q \in N(v)$, $y_Q \in V(Q)$ and $e_Q \in A_Q$. Let $A_{\mathcal{H}} = \{e_Q \colon Q \in \mathcal{H}\}$. For each $Q \in \mathcal{H}$, we have $e_Q \cap N(v) \neq \emptyset$ and $|e_Q| = r$, so $e_Q$ links at most $r-1$ distinct members of $\mathcal{H}$ to $N(v)$. Thus, $|\mathcal{H}| \leq (r-1)|A_{\mathcal{H}}|$, and hence $|A_{\mathcal{H}}| \geq \frac{\mathcal{|H}|}{r-1}$. For any hypergraph $Q$, let $m_Q = |E(Q)|$, and let $D_Q$ be a $K_k^r$-isolating set of $Q$ of size $\iota(Q,K^r_k)$. Clearly, $A_v$, $A_{\mathcal{H}}$ and $\bigcup_{Q \in \mathcal{H}} E(Q)$ are pairwise disjoint subsets of $E(H)$, so by (\ref{eq1}),
\begin{equation} m \ge k_r + 1 + \frac{\mathcal{|H}|}{r-1} + \sum \limits_{Q\in\mathcal{H}} m_Q. \label{eq3}     
\end{equation}
Let $\mathcal{H'} = {\rm C}_k^r(H')$ (so $\mathcal{H}' = \{Q \in \mathcal{H} \colon Q \simeq K^r_k\}$). By the induction hypothesis, $|D_Q| \le ((r-1)m_Q + 1)/t$ for each $Q \in \mathcal{H} \setminus \mathcal{H}'$.

Suppose $\mathcal{H'} = \emptyset$. By Lemma~\ref{lemma1} (with $X= \{v\}$, $Y = N[v]$ and $Z = \emptyset$), Lemma~\ref{lemma2} and (\ref{eq3}), we have
\begin{align} t \left( \iota(H,K^r_k) \right) & \le t \left( 1 + \iota(H',K^r_k) \right) = t \bigg{(} 1 + \sum\limits_{Q \in \mathcal{H}} |D_Q| \bigg{)} \leq t +  \sum_{Q \in \mathcal{H}} ((r-1)m_Q + 1) \nonumber\\
&= (r-1) \bigg{(} k_r + 1 + \bigg{(} \sum_{Q \in \mathcal{H}} m_Q \bigg{)} + \frac{|\mathcal{H}|}{r-1} \bigg{)} + 1 \leq (r-1)m + 1, \nonumber %\label{case1}
\end{align}
so $\iota(H,K^r_k) \leq ((r-1)m + 1)/t$.

Now suppose $\mathcal{H'} \neq \emptyset$. Let $Q_1 \in \mathcal{H'}$. Let $e_1 = e_{Q_1}$, $x = x_{Q_1}$ and $y = y_{Q_1}$. %Removing the edge $e_1$ to obtain 
Let $H_1 = H - (\{e_1\} \cup E(Q_1))$. Then, $H_1$ has a component $H_v^*$ with $N[v] \subseteq V(H_v^*)$, and each member of ${\rm C}(H_1) \setminus \{H_v^*\}$ is a member of $\mathcal{H}$ that is linked to $N(v)$ only by $e_1$, or a $1$-vertex subhypergraph of $Q_1$. Let $\mathcal{H}_1 = {\rm C}(H_1) \cap \mathcal{H}'$, $\mathcal{H}_2 = {\rm C}(H_1) \cap (\mathcal{H} \setminus \mathcal{H}')$ and $\mathcal{H}_{Q_1} = \{Z \in {\rm C}(H_1) \colon V(Z) \in {V(Q_1) \choose 1}\}$. Then, ${\rm C}(H_1) = \{H_v^*\} \cup \mathcal{H}_1 \cup \mathcal{H}_2 \cup \mathcal{H}_{Q_1}$. Let $h_1 = |\mathcal{H}_1|$ and $h_2 = |\mathcal{H}_2|$. We have $|e_1| = r$, $e_1 \cap N(v) \neq \emptyset$, $e_1 \cap V(Q) \neq \emptyset$ for each $Q \in \{Q_1\} \cup \mathcal{H}_1 \cup \mathcal{H}_2$, and $Q_1 \notin \mathcal{H}_1 \cup \mathcal{H}_2$, so $h_1 + h_2 \leq r - 2$. 

Let $\mathcal{H}_2^* = \{H_v^*\} \cup \mathcal{H}_2$. Since $E(H) = \{e_1\} \cup E(Q_1) \cup E(H_1)$,
\begin{equation} m = k_r + 1 + h_1 k_r + \sum_{Q \in \mathcal{H}_2^*} m_Q. \label{eq4}
\end{equation}
By (\ref{eq1}), $H_v^* \not\simeq K_k^r$. By the induction hypothesis, $\iota(H_v^*,K^r_k) \le ((r-1)m_{H_v^*} + 1)/t$. Let $D = \{y\} \cup \bigcup_{Q \in \mathcal{H}_2^*} D_Q$. Since $r \leq k$, for each $z \in V(Q_1)$, we have $y, z \in e_{z}$ for some $e_{z} \in E(Q_1)$, so $z \in N[y]$. Thus, $V(Q_1) \subseteq N[y]$, and hence $N[y] \cap e \neq \emptyset$ for each $e \in E(Q_1)$. Since $y \in e_1$, we have $e_1 \subseteq N[y]$, and hence $N[y] \cap V(Q) \neq \emptyset$ for each $Q \in \mathcal{H}_1$. Thus, $D$ is a $K^r_k$-isolating set of $H$, and hence 
\begin{align}
t \left( \iota(H,K^r_k) \right) &\leq t|D| = t + \sum_{Q \in \mathcal{H}_2^*} t |D_Q| \leq t + \sum_{Q \in \mathcal{H}_2^*} ((r-1)m_Q + 1) \nonumber\\
&= (r-1) \bigg{(} k_r + 1 + \sum_{Q \in \mathcal{H}_2^*}m_Q \bigg{)} + |\mathcal{H}_2^*| + 1 \nonumber \\
&= (r-1)(m - h_1 k_r) + h_2 + 2. \label{4}
\end{align}
\medskip

\noindent \textbf{Case~1:} $h_1 \ge 1$. Since $h_1 + h_2 \leq r-2$, $h_2 \le r-3$. Since $r \leq k$, $k_r \ge 1$. By (\ref{4}),
\begin{equation} t \left( \iota(H,K^r_k) \right) \leq (r-1)m - (r-1)h_1 k_r + r - 1 = (r-1)m + (r-1)(1 - h_1 k_r) \leq (r-1)m, \nonumber
\end{equation}
so $\iota(H,K^r_k) \leq (r-1)m/t$.
\\

\noindent \textbf{Case~2:} $h_1 = 0$. For any $Q \in \mathcal{H}_2$ and $w \in e_1 \cap V(Q)$, let $A_w = E_Q(w)$. Recall that $e_1$ is the only member of $E(H)$ that links $Q$ to $N[v]$, so 
\begin{equation} E_H(w) = \{e_1\} \cup A_w. \nonumber
\end{equation}
Since $r \leq k$, $y \in e_y$ for some $e_y \in E(Q_1)$. Let $Q_1' = Q_1 - \{e_y\}$. Since $Q_1 \simeq K_k^r$, $Q_1'$ contains no $K_k^r$-copy.
\\

\noindent\textbf{Case~2.1:} For some $Q_2 \in \mathcal{H}_2$ and some $w \in e_1 \cap V(Q_2)$, $Q_2 - A_{w}$ has a component $C$ that is a $K^r_k$-copy. Since $Q_2$ is connected, we have $w, w_C \in e_C$ for some $w_C \in V(C)$ and $e_C \in A_w$. Let $Q_2' = Q_2 - A_w$, $W = (\{w\}, \emptyset)$, $\mathcal{H}_3 = {\rm C}_k^r(Q_2')$ and $\mathcal{H}_4 = {\rm C}(Q_2') \setminus (\{W\} \cup \mathcal{H}_3)$ (note that $W \in {\rm C}(Q_2')$). Let $H_2 = H - (\{e_1, e_y\} \cup A_w)$. Then, $H_2$ has a component $H_v'$ with $N[v] \subseteq V(H_v')$, and ${\rm C}(H_2) = \{H_v', W\} \cup (\mathcal{H}_2 \setminus \{Q_2\}) \cup \mathcal{H}_3 \cup \mathcal{H}_4 \cup \mathcal{H}_{Q_1'}$, where $\mathcal{H}_{Q_1'}$ is a set (possibly empty) of subhypergraphs of $Q_1'$. For any $e \in A_w$ and $I \in {\rm C}(Q_2')$, we say that \emph{$e$ links $I$ to $w$} if $e \cap V(I) \neq \emptyset$. Since $Q_2$ is connected, for each $I \in {\rm C}(Q_2')$, some $e_I \in A_w$ links $I$ to $w$. Since $w, w_C \in e_C$, $e_C$ links at most $r-2$ distinct members of $\mathcal{H}_4$ to $w$. For each $e \in A_w \setminus \{e_C\}$, we have $w \in e$, so $e$ links at most $r-1$ distinct members of $\mathcal{H}_4$ to $w$. Thus, $|\mathcal{H}_4| \leq (r-2) + (r-1)(|A_w| - 1) = (r-1)|A_w| - 1$, and hence $|A_w| \geq \frac{|\mathcal{H}_4| + 1}{r-1}$. Let $\mathcal{I} = \{H_v'\} \cup (\mathcal{H}_2 \setminus \{Q_2\}) \cup \mathcal{H}_4$. Thus, 
\begin{equation} |\mathcal{I}| = h_2 + |\mathcal{H}_4| \leq r - 2 + |\mathcal{H}_4|. \label{eq:4.5}
\end{equation}
%sine |H_2\Q|\leq h2-1 and but +1 for H'_v. 
Since $E(H) = \{e_1, e_y\} \cup A_w \cup E(H_2) \supseteq \{e_1, e_y\} \cup A_w \cup \bigcup_{I \in \mathcal{I} \cup \mathcal{H}_3} E(I)$ and $C \in \mathcal{H}_3$,
\begin{equation} m \ge k_r + 2 + \frac{|\mathcal{H}_4| + 1}{r-1} + \sum_{I \in \mathcal{I}} m_I. \label{eq:5}
\end{equation}
By (\ref{eq1}), $H_v' \not\simeq K_k^r$. By the induction hypothesis, $|D_I| \leq ((r-1)m_I + 1)/t$ for each $I \in \mathcal{I}$. Let $D = \{w\} \cup \bigcup_{I \in \mathcal{I}} D_I$. Since $y, w \in e_1$, we have $y \in N[w]$, so $N[w] \cap e_y \neq \emptyset$. For each $I \in \mathcal{H}_3$, $\emptyset \neq e_I \cap V(I) \subseteq N[w] \cap V(I)$. Thus, $D$ is a $K^r_k$-isolating set of $H$. %since $N[x_Q]$ eliminates any $k$-clique disconnected by the removal of $x_Q$. 
Therefore,
\begin{align} t \left( \iota(H,K^r_k) \right) &\leq t|D| = t + \sum_{I \in\mathcal{I}} t|D_I| \leq t + \sum_{I \in \mathcal{I}} ((r-1)m_I + 1) \nonumber \\
&= (r-1) \bigg{(} k_r + 1 + \sum_{I \in \mathcal{I}} m_I \bigg{)} + 1 + |\mathcal{I}| \nonumber \\
&\le (r-1) \bigg{(} m - 1 - \frac{|\mathcal{H}_4| + 1}{r-1} \bigg{)} + r-1 + |\mathcal{H}_4| \quad \mbox{(by (\ref{eq:4.5}) and (\ref{eq:5}))} \nonumber \\
&= (r-1)m - 1, \nonumber
\end{align}
and hence $\iota(H, K_k^r) \leq ((r-1)m - 1)/t$.
\\

\noindent\textbf{Case~2.2:} For every $Q \in \mathcal{H}_2$ and $w \in e_1 \cap V(Q)$, no component of $Q - A_{w}$ is a $K_k^r$-copy. For each $Q\in \mathcal{H}_2$, let $w_Q \in e_1 \cap V(Q)$, $Q' = Q - A_{w_Q}$, $W_Q = (\{w_Q\}, \emptyset)$ and $\mathcal{I}_Q = {\rm C}(Q') \setminus \{W_Q\}$ (note that $W_Q \in {\rm C}(Q')$). Let $\mathcal{H}_3 = \bigcup_{Q \in \mathcal{H}_2} \mathcal{I}_Q$. Consider any $Q \in \mathcal{H}_2$ and any $I \in \mathcal{I}_Q$. Then, $I$ is not a $K_k^r$-copy. For any $e \in A_{w_Q}$, we say that \emph{$e$ links $I$ to $w_Q$} if $e\cap V(I) \neq \emptyset$. Since $Q$ is connected, some $e_I \in A_{w_Q}$ links $I$ to $w_Q$. For each $e \in A_{w_Q}$, we have $w_Q \in e$, so $e$ links at most $r-1$ members of $\mathcal{I}_Q$ to $w_Q$. Thus, $|\mathcal{H}_3| \leq \sum_{Q\in \mathcal{H}_2} (r-1) |A_{w_Q}|$, and hence 
\begin{equation} \sum_{Q \in \mathcal{H}_2} |A_{w_Q}| \ge \frac{|\mathcal{H}_3|}{r-1}. \label{ewqsum}
\end{equation}
\medskip

\noindent \textbf{Case~2.2.1:} For some $e_2 \in E(H) \setminus \{e_1\}$, $e_2$ links $Q_1$ to $N(v)$. Thus, $y_2 \in e_2$ for some $y_2 \in V(Q_1)$. Let 
\[\mathcal{H}_4 = \{Q \in \mathcal{H} \colon e_2 \mbox{ is the only member of $E(H)$ that links $Q$ to } N[v]\}.\]
Note that $(\{Q_1\} \cup \mathcal{H}_2) \cap \mathcal{H}_4 = \emptyset$. Let $M = \{e_1,e_2\} \cup E(Q_1) \cup \bigcup_{Q \in \mathcal{H}_2} A_{w_Q}$ and $H_2 = H - M$. Then, $H_2$ has a component $H'_v$ with $N[v]\subseteq V(H_v')$, and ${\rm C}(H_2) = \{H_v'\} \cup \mathcal{H}_3 \cup \mathcal{H}_4 \cup \{W_Q \colon Q \in \mathcal{H}_2\} \cup \mathcal{H}_{Q_1}'$, where $\mathcal{H}_{Q_1}'$ is a set (possibly empty) of $1$-vertex subhypergraphs of $Q_1$. 

If some member of $\mathcal{H}_4$ is a $K_k^r$-copy, then the result follows by the argument in Case~1 with $e_2$ and $\mathcal{H}_4$ now taking the role of $e_1$ and $\mathcal{H}_1 \cup \mathcal{H}_2$, respectively, in Case~1. Suppose that no member of $\mathcal{H}_4$ is a $K_k^r$-copy. We have $|e_2| = r$, $y_2 \in e_2$ and $e_2 \cap N[v] \neq \emptyset$, so $|\mathcal{H}_4| \leq r-2$.  

Let $\mathcal{I} = \{H_v'\} \cup \mathcal{H}_3 \cup \mathcal{H}_4$. Then, 
\begin{equation}
    |\mathcal{I}|= 1+ |\mathcal{H}_3| + |\mathcal{H}_4|\leq r-1 + |\mathcal{H}_3|. \label{eq:4.7}
\end{equation}
We have $E(H) = M \cup E(H_2) = M \cup \bigcup_{I \in \mathcal{I}} E(I)$, so by (\ref{ewqsum}),
\begin{equation} m \ge k_r + 2 + \frac{|\mathcal{H}_3|}{r-1} + \sum_{I \in \mathcal{I}} m_I. \label{eq:4.8}
\end{equation}
By (\ref{eq1}), $H_v' \not\simeq K_k^r$. By the induction hypothesis, $|D_I| \leq (m_I(r-1)+1)/t$ for each $I \in \mathcal{I}$. Let $D = \{y\} \cup \bigcup_{I \in \mathcal{I}} D_I$. Recall that $N[y] \cap e \neq \emptyset$ for each $e \in E(Q_1)$. Since $V(Q_1) \subseteq N[y]$, $y_2 \in N[y] \cap e_2$. Since $\{w_Q \colon Q \in \mathcal{H}_2\} \subset e_1 \subseteq N[y]$, $N[y] \cap e \neq\emptyset$ for each $e \in \bigcup_{Q \in \mathcal{H}_2} A_{w_Q}$. By Lemma~\ref{lemma1} (with $X = \{y\}$, $Y = \emptyset$, $Z = M$ and $S = \bigcup_{I \in \mathcal{I}} D_I$) and Lemma~\ref{lemma2} (applied to $H_2$), $D$ is a $K_k^r$-isolating set of $H$. Therefore,
\begin{align}
    t(\iota(H,K^r_k))\leq t|D| &= t + \sum_{I\in\mathcal{I}} t|D_I| \leq t +\sum_{I \in \mathcal{I}}((r-1)m_I + 1) \nonumber \\
    &=(r-1) \bigg{(} k_r + 1 + \sum_{I \in \mathcal{I}} m_I \bigg{)} + 1 + |\mathcal{I}| \nonumber \\
    &\leq (r-1) \bigg{(} m - 1 - \frac{|\mathcal{H}_3|}{r-1} \bigg{)} + r + |\mathcal{H}_3| \quad \mbox{(by (\ref{eq:4.7}) and (\ref{eq:4.8}))} \nonumber \\
    &= (r-1)m + 1, \nonumber 
\end{align}
and hence $\iota(H, K_k^r) \leq ((r-1)m + 1)/t$.\\

\noindent\textbf{Case~2.2.2:} No member of $E(H) \setminus \{e_1\}$ links $Q_1$ to $N(v)$. Then, since $E(H-\{e_1\}) = E(H_1) \cup E(Q_1)$ and $h_1 = 0$, ${\rm C}(H - \{e_1\}) = \{H_v^*, Q_1\} \cup \mathcal{H}_2$. Let $M = \{e_1,e_y\}  \cup \bigcup_{Q \in \mathcal{H}_2} A_{w_Q}$ and $H_2 = H - M$. Then, ${\rm C}(H_2) = \{H_v^*\} \cup {\rm C(Q_1')} \cup \{W_Q \colon Q \in \mathcal{H}_2\} \cup \mathcal{H}_3$. 

Recall that $x \in e_1 \cap N(v)$. Let $A_x = E_{H_v^*}(x)$, $H_v' = H_v^* - A_x$, $W_x = (\{x\}, \emptyset)$, $\mathcal{H}_4 = {\rm C}^r_k(H_v')$ and $\mathcal{H}_5 = {\rm C}(H_v') \setminus (\{W_x\} \cup \mathcal{H}_4)$ (note that $W_x \in {\rm C}(H_v')$). For any $e \in A_x$ and $I \in \mathcal{H}_4 \cup \mathcal{H}_5$, we say that \emph{$e$ links $I$ to $x$} if $e \cap V(I) \neq \emptyset$. Since $H_v^*$ is connected, for each $I \in\mathcal{H}_4 \cup \mathcal{H}_5$, some $e_I \in A_x $ links $I$ to $x$. For each $e \in A_x$, we have $x \in e$, so $e$ links at most $r-1$ members of $\mathcal{H}_5$ to $x$. Thus, $|\mathcal{H}_5| \leq (r-1)|A_x|$. Let $\mathcal{I} = \mathcal{H}_3 \cup \mathcal{H}_5$. We have $E(H) = M \cup E(H_2) \supseteq E(Q_1) \cup \{e_1\} \cup \left(\bigcup_{Q \in \mathcal{H}_2} A_{w_Q} \right) \cup A_x \cup \bigcup_{I\in \mathcal{I}\cup \mathcal{H}_4 } E(I)$. Together with (\ref{ewqsum}), this yields
\begin{equation} m \ge k_r + 1 + \frac{|\mathcal{H}_3|}{r-1} + \frac{|\mathcal{H}_5|}{r-1} + \sum_{I \in \mathcal{I}} m_I = k_r + 1 + \frac{|\mathcal{I}|}{r-1} + \sum_{I \in \mathcal{I}} m_I. \label{eq:4.10}
\end{equation}
By the induction hypothesis, $|D_I| \leq ((r-1)m_I + 1)/t$ for each $I \in \mathcal{I}$. Let $D_v^* = \{x\} \cup \bigcup_{I \in \mathcal{H}_5} D_I$ and $D = D_v^* \cup \bigcup_{I \in \mathcal{H}_3} D_I$. For each $I \in \mathcal{H}_4$, $\emptyset \neq e_I \cap V(I) \subseteq N_{H_v^*}[x] \cap V(I)$. Thus, $D_v^*$ is a $K_k^r$-isolating set of $H_v^*$, and hence $D$ is a $K_k^r$-isolating set of $H_2$ (by Lemma~\ref{lemma2}). Since $x \in e_1$, $e_1 \subseteq N[x]$. Since $y \in e_1 \cap e_y$, $y \in N[x] \cap e_y$. Since $\{w_Q \colon Q \in \mathcal{H}_2\} \subset e_1 \subseteq N[x]$, $N[x] \cap e \neq \emptyset$ for each $e \in \bigcup_{Q \in \mathcal{H}_2} A_{w_Q}$. By Lemma~\ref{lemma1} (with $X = \{x\}$, $Y = \emptyset$, $Z = M$ and $S = D$), $D$ is a $K_k^r$-isolating set of $H$. Therefore,
\begin{align}
    t(\iota(H,K^r_k)) &\leq t|D| = t  + \sum_{I \in \mathcal{I}} t |D_I| \leq t + \sum_{I \in \mathcal{I}} ((r-1)m_I + 1) \nonumber \\
    &=(r-1) \bigg{(} k_r + 1 + \sum_{I \in \mathcal{I}} m_I \bigg{)} + 1 + |\mathcal{I}| \nonumber \\
    &\le (r-1) \bigg{(} m - \frac{|\mathcal{I}|}{r-1} \bigg{)} + 1 + |\mathcal{I}| \quad \mbox{(by (\ref{eq:4.10}))} \notag\\
    &= (r-1)m + 1, \nonumber
\end{align}
and hence $\iota(H, K_k^r) \leq ((r-1)m + 1)/t$.~\hfill{$\Box$}
\\
 
\noindent \textbf{Proof 2 of Theorem~\ref{edgethm} for $r < k$}. As shown in Proof~1, the bound is attained if $H = B_{m,k,r}$. We now provide an alternative proof that the bound holds for $r < k$, again using induction on $m$. If $r = 2$, then the result is given by (\ref{BFK_bound_m}). Suppose $r \ge 3$. Then, $k \ge 4$. Let $k_r = \binom{k}{r}$ and $t = (r-1)(k_r+1) + 1$. The result is trivial if $m \leq 2$ or $H$ has no $k$-cliques. Suppose that $m \ge 3$ and $H$ has a $k$-clique $C$. Since $H$ is connected and $H \not\simeq K_k^r$, $N[v] \setminus V(C) \neq \emptyset$ for some $v \in V(C)$. Let $u \in N[v] \setminus V(C)$, and let $e_u \in E(H)$ such that $v, u \in e_u$. Let $W = V(C)\cup e_u$ and $A_W = E(H[W])$. Thus, $E(C)\cup \{e_u\} \subseteq A_W$. We have
\begin{equation} |W| \ge k+1  \quad \text{and} \quad |A_W| \geq k_r + 1. \label{eq:2.1}
\end{equation}
If $V(H) = W$, then $\{v\}$ is a $k$-clique isolating set of $H$, so $\iota(H,K^r_k) = 1 \leq ((r-1)m + 1)/t$ by (\ref{eq:2.1}).
Suppose $V(H)\neq W$. Let $H' = H - W$. Then, $V(H') \neq \emptyset$. Let $\mathcal{H} = {\rm C}(H')$. For each $Q \in \mathcal{H}$, let $A_Q = \{e \in E(H) \colon e \cap (W) \neq \emptyset \neq e \cap V(Q)\}$. 

Suppose $e \notin A_W \cup \bigcup_{Q \in \mathcal{H}} (E(Q) \cup A_Q)$ for some $e \in E(H)$. Let $\mathcal{H}_e = \{Q \in \mathcal{H} \colon e \cap V(Q) \neq \emptyset\}$. Since $e \notin A_W$, we have $e \setminus W \neq \emptyset$, so $ \mathcal{H}_e \neq \emptyset$. Let $Q_e \in \mathcal{H}_e$. Since $e \notin A_{Q_e}$, $e \cap W = \emptyset$. Thus, $e \subseteq V(H')$. Let $H'_e = (\bigcup_{Q \in \mathcal{H}_e} V(Q), \{e\} \cup \bigcup_{Q \in \mathcal{H}_e} E(Q))$. Then, $H'_e$ is a connected subhypergraph of $H'$. Since $e \notin E(Q_e)$, we have $E(Q_e) \subsetneq E(H'_e)$, which contradicts $Q_e \in \mathcal{H}$. Therefore,
\begin{equation} E(H) = A_W \cup \bigcup_{Q \in \mathcal{H}} (E(Q) \cup A_Q). \label{H_edges}
\end{equation}

For any $Q \in \mathcal{H}$ and $x \in W$, we say that $Q$ is \emph{linked to $x$} if $x \in e$ for some $e \in A_Q$, and we say that $e$ \emph{links $Q$ to $x$} or that $e$ \emph{links $Q$ to $W$}. Since $H$ is connected, each member of $\mathcal{H}$ is linked to at least one member of $W$. Thus, for each $Q \in \mathcal{H}$, we have $x_Q, y_Q \in e_Q$ for some $x_Q \in W$, $y_Q \in V(Q)$ and $e_Q \in A_Q$. Let $A_{\mathcal{H}} = \{e_Q \colon Q \in \mathcal{H}\}$. For each $Q \in \mathcal{H}$, we have $e_Q \cap W\neq \emptyset$ and $|e_Q| = r$, so $e_Q$ links at most $r-1$ distinct members of $\mathcal{H}$ to $W$. Thus, $|\mathcal{H}| \leq (r-1)|A_{\mathcal{H}}|$, and hence $|A_{\mathcal{H}}| \geq \frac{\mathcal{|H}|}{r-1}$. For any hypergraph $Q$, let $m_Q = |E(Q)|$, and let $D_Q$ be a $K_k^r$-isolating set of $Q$ of size $\iota(Q,K^r_k)$. Clearly, $A_W$, $A_{\mathcal{H}}$ and $\bigcup_{Q \in \mathcal{H}} E(Q)$ are pairwise disjoint subsets of $E(H)$, so by (\ref{eq:2.1}),
\begin{equation} |E(H)| \ge k_r + 1 + \frac{\mathcal{|H}|}{r-1} + \sum \limits_{Q\in\mathcal{H}} m_Q . \label{eq2.4}     
\end{equation}
Let $\mathcal{H'} = {\rm C}_k^r(H')$. By the induction hypothesis, $|D_Q| \le ((r-1)m_Q+1)/t$ for each $Q \in \mathcal{H} \setminus \mathcal{H}'$. If $\mathcal{H'} = \emptyset$, then we obtain $\iota(H,K^r_k) \leq ((r-1)m + 1)/t$ as in Proof 1 of Theorem~\ref{edgethm}. Suppose $\mathcal{H'}\neq\emptyset$. Let $Q_1\in \mathcal{H'}$. Let $e_1 = e_{Q_1}$, $x = x_{Q_1}$ and $y = y_{Q_1}$. Let $e_C \in E(C)$ such that 
$x \in e_C$ if $x \in V(C)$, and $v \in e_C$ if $x \notin V(C)$. Let $C' = C - \{e_C\}$ and 
$$H_1 = H - (\{e_1, e_C\} \cup E(Q_1)).$$

Consider any $z_1, z_2 \in V(C)$ with $z_1 \neq z_2$. If $\{z_1, z_2\} \nsubseteq e_C$, then $C'$ has ${k-2 \choose r-2}$ edges containing both $z_1$ and $z_2$. Suppose $\{z_1, z_2\} \subseteq e_C$. Since $r \geq 3$, $e_C \setminus \{z_1, z_2\} \neq \emptyset$. Let $z_3 \in e_C \setminus \{z_1, z_2\}$. Since $r < k$, $V(C) \setminus e_C \neq \emptyset$. Let $z_4 \in V(C) \setminus e_C$. We have $z_1, z_2 \in (e_C \setminus \{z_3\}) \cup \{z_4\} \in E(C')$. Thus, $C'$ is connected. 

Therefore, $H_1$ has a component $H_{W}^*$ with $W \subseteq V(H_W^*)$, and each member of ${\rm C}(H_1) \setminus \{H_W^*\}$ is a member of $\mathcal{H}$ that is linked to $W$ only by $e_1$, or a $1$-vertex subhypergraph of $Q_1$. Let $\mathcal{H}_1=\ {\rm C}(H_1) \cap \mathcal{H'}$, $\mathcal{H}_2 = {\rm C}(H_1) \cap (\mathcal{H} \setminus \mathcal{H'})$ and $\mathcal{H}_{Q_1} = \{Z \in {\rm C}(H_1) \colon V(Z) \in {V(Q_1) \choose 1}\}$. Then, 
\begin{equation} {\rm C}(H_1) = \{H_W^*\} \cup \mathcal{H}_1 \cup \mathcal{H}_2 \cup \mathcal{H}_{Q_1}. \label{H_1comps}
\end{equation}
Let $h_1 = |\mathcal{H}_1|$ and $h_2 = |\mathcal{H}_2|$. We have $|e_1| = r$, $e_1 \cap W \neq \emptyset$, $e_1 \cap V(Q) \neq \emptyset$ for each $Q \in \{Q_1\} \cup \mathcal{H}_1 \cup \mathcal{H}_2$, and $Q_1 \notin \mathcal{H}_1 \cup \mathcal{H}_2$, so $h_1 + h_2 \leq r - 2$. Let $\mathcal{H}_2^* = \{H_W^*\} \cup \mathcal{H}_2$. Then, 
\begin{equation} |\mathcal{H}_2^*| \leq r-1. \label{H_2^*ineq}
\end{equation}
Since $E(H) = E(Q_1) \cup \{e_1,e_C\} \cup E({H}_1) $, 
\begin{equation} m = k_r+2 + h_1 k_r + \sum\limits_{Q \in \mathcal{H}_2^*} m_{Q}. \label{eq:2.3}
\end{equation}
By (\ref{eq:2.1}), $|V(H_W^*)| \ge k+1$, so $H_W^* \not\simeq K_k^r$. By the induction hypothesis, $\iota(H_W^*,K_k^r) \leq ((r-1)m_{H_W^*} + 1)/t$. \\

\noindent\textbf{Case~1:} $x \in V(C)$. Thus, $x \in e_C$. Let $D = \{y\} \cup \bigcup_{Q\in \mathcal{H}_2^*} D_Q$. We have $x \in N[y] \cap e_C$. Since $r < k$, for each $z\in V(Q_1)$, we have $y, z\in e_z$ for some $e_z\in E(Q_1)$, so $z\in N[y]$. Thus, $V(Q_1)\subseteq N[y]$, and hence $N[y] \cap e \neq \emptyset$ for each $e\in E(Q_1)$. Since $y \in e_1$, we have $e_1 \subseteq N[y]$, so $N[y] \cap V(Q) \neq \emptyset$ for each $Q\in \mathcal{H}_1$. Thus, $D$ is a $K^r_k$-isolating set of $H$. Therefore, 
\begin{align} t(\iota(H,K^r_k)) &\leq t|D| = t + \sum_{Q \in \mathcal{H}^*_2} t|D_Q| \leq (r-1)\bigg{(} k_r + 1 + \sum_{Q \in \mathcal{H}_2^*} m_Q \bigg{)} + 1 + |\mathcal{H}_2^*| \nonumber \\
&\leq (r-1)(m - 1 - h_1 k_r) + r \quad \mbox{(by (\ref{H_2^*ineq}) and (\ref{eq:2.3})}) \nonumber \\
&\leq (r-1)(m-1) + r = (r-1)m + 1, \nonumber
\end{align}
and hence $\iota(H,K^r_k) \leq ((r-1)m + 1)/t$.
\\

\noindent\textbf{Case~2:} $x \notin V(C)$. Thus, $v \in e_C$ and $x \in e_u$. If $e \cap V(C) \neq \emptyset$ for some $e \in A_{Q_1}$, then the result follows by the argument in Case~1 with $e$ and any $x' \in e \cap V(C)$ and $y' \in e \cap V(Q_1)$ taking the role of $e_1$, $x$ and $y$, respectively. Now suppose $e \cap V(C) = \emptyset$ for each $e \in A_{Q_1}$. Then, 
\begin{equation} e \cap e_u \neq \emptyset \mbox{ for each } e \in A_{Q_1}. \label{e_u_int}
\end{equation}

Let $D = \{x\} \cup \bigcup_{Q \in \mathcal{H}_2^*} D_Q$. As in Case~1, $t|D| \leq (r-1)m + 1$, so the result follows if $D$ is a $K^r_k$-isolating set of $H$. Suppose that $D$ is not a $K^r_k$-isolating set of $H$. Then, $N[D] \cap V(C^*) = \emptyset$ for some $C^* \in \mathcal{C}_k^r(H)$. We have $E(H) = E(H_1) \cup \{e_1, e_C\} \cup E(Q_1)$. Since $x \in e_1$, we have $e_1 \subseteq N[x]$, and hence $N[x] \cap V(Q) \neq \emptyset$ for every $Q \in \{Q_1\} \cup \mathcal{H}_1$. We have $v \in e_u \cap e_C \subseteq N[x] \cap e_C$. Suppose $E(C^*) \cap E(Q_1) = \emptyset$. Then, $E(C^*) \subseteq E(H_1)$ and $C^* \notin \mathcal{H}_1$. By (\ref{H_1comps}), $C^* \in \mathcal{C}_k^r(Q)$ for some $Q \in \mathcal{H}_2^*$. We have $N[D_Q] \cap V(C^*) \neq \emptyset$, which contradicts $N[D] \cap V(C^*) = \emptyset$. Thus, $E(C^*) \cap E(Q_1) \neq \emptyset$. Since $\emptyset \neq N[x] \cap V(Q_1) \subseteq N[D] \cap V(Q_1)$, $Q_1 \neq C^*$. Thus, since $C^* \simeq K_k^r$, $e^* \cap V(Q_1) \neq \emptyset$ for some $e^* \in E(C^*) \setminus E(Q_1)$. By (\ref{H_edges}), $e^* \in A_{Q_1}$. By (\ref{e_u_int}), we have $\emptyset \neq e_u \cap e^* \subseteq N[x] \cap e^*$, which contradicts $N[D] \cap V(C^*) = \emptyset$.~\hfill{$\Box$}

\end{document}